\documentclass[a4paper,12pt]{article}

\usepackage[table,dvipsnames,svgnames,x11names,hyperref]{xcolor}
\usepackage[utf8]{inputenc}

\usepackage{lmodern}
\usepackage[T1]{fontenc}

\usepackage[top=1in, bottom=1.25in, left=1.25in, right=1.25in]{geometry}
\usepackage{amsmath,amssymb,amsthm}

\usepackage{tikz}
\usepackage{array}
\usepackage{enumerate}
\usepackage{graphicx}
\usepackage{float}
\usepackage{caption}
\usepackage{subcaption}
\usepackage[colorlinks=true, allcolors=MidnightBlue]{hyperref}
\usepackage{tasks}

\newtheorem{thm}{Theorem}
\newtheorem{cor}[thm]{Corollary}
\newtheorem{lem}[thm]{Lemma}
\newtheorem{prop}[thm]{Proposition}

\theoremstyle{definition}

\theoremstyle{remark}

\newtheorem*{rem}{Remark}

\newtheorem*{example}{Example}

\author{Sel\c{c}uk Kayacan\thanks{I would like to thank to Volkmar Welker for the helpful suggestions he made during the preparation of this paper.}}
\title{Some remarks on the subrack lattice of finite racks}
\date{}

\begin{document}

\maketitle

\small

\begin{center}
  Bah\c{c}e\c{s}ehir University, Faculty of Engineering\\ and Natural Sciences,
  Istanbul, Turkey\\
  {\it e-mail:} \href{mailto:selcuk.kayacan@eng.bau.edu.tr}{selcuk.kayacan@eng.bau.edu.tr}
\end{center}

\begin{abstract}
       The set of all subracks $\mathcal{R}(X)$ of a finite rack $X$ form a lattice under inclusion. We prove that if a rack $X$ satisfies a certain condition then the homotopy type of the order complex of $\mathcal{R}(X)$ is a $(m-2)$-sphere, where $m$ is the number of maximal subracks of $X$. The rack $X$ satisfying the condition of this general result is necessarily decomposable. Two particular instances occur when
       \begin{itemize}
       \item $X=G$ is a group rack, and when
       \item $X=C$ is a conjugacy class rack of a nilpotent group.
       \end{itemize}
We also studied the subrack lattices of indecomposable racks by focusing on the conjugacy class racks of symmetric or alternating groups and determined the homotopy types of the corresponding order complexes in some cases.


\end{abstract}

\section{Introduction}\label{sec:intro}

The study of combinatorial objects associated to algebraic structures is a recurring theme in literature. Subgroup lattices is among the most prominent research areas in this regard \cite{Sch94}. One may ask, for example, which classes of groups can be distinguished by their subgroup lattices in this combinatorial setting. An alternative approach is adopted in the pioneering works of Brown \cite{Bro75} and Quillen \cite{Qui78}. Those authors obtained striking results by introducing the order complexes of some specific subgroup posets (partially ordered sets) which renders the use of topological methods possible. Many other works were published subsequently in this new setting of subgroup complexes \cite{Smi11}.

Still some variations on this general framework is possible. As an example we may consider the intersection graph of subgroups of a group which is intimately related with the subgroup lattice of the group. Actually, given the subgroup lattice one can recover the intersection graph but not vice-versa. In \cite{Kay18} the intersection complex of a group is introduced as the simplicial complex whose faces are the sets of proper subgroups which intersect non-trivially. The intersection complex share the same homotopy type with the order complex of the subgroup lattice and its 1-skeleton is the intersection graph. Moreover, it is contractible if domination number of the intersection graph is 1. 

In \cite{HSW19} Heckenberger, Shareshian, and Welker initiated the study of racks from the combined perspective of combinatorics and finite group theory. A \emph{rack} $X$ is a set (possibly empty) together with a binary operation $\triangleright$ satisfying the following two axioms:
\begin{itemize}
\item[\textbf{(A1)}] for all $a,b,c\in X$ we have $a \triangleright (b \triangleright c) = (a \triangleright b) \triangleright (a \triangleright c)$ 
\item[\textbf{(A2)}] for all $a,b\in X$ there is a
  unique $x\in X$ such that $a \triangleright x = b$.
\end{itemize}
Since (A2) tells us $\phi_a\colon X\to X;\,x\mapsto a\triangleright x$ is a bijection for every $a\in X$, and since (A1) (self distributivity) means that $\phi_a$ respects to the structure of $X$, an alternative definition of a rack would be that it is a set with a binary operation such that left multiplications are automorphisms of the set. We say $X$ is a \emph{quandle} if, in addition, it satisfies the following third axiom:
\begin{itemize}
\item[\textbf{(A3)}] for every $a\in X$ we have $a\triangleright a = a$.
\end{itemize}
For example, any finite group $G$ or any conjugacy class $C$ of $G$ equipped with the conjugation operation $a\triangleright b := aba^{-1}$ satisfies all the above three axioms; hence, is an example of a \emph{quandle}. In the first case, we say $G$ is a \emph{group rack} and in the latter case we say $C$ is a \emph{conjugacy class rack}. More generally, we say $X$ is a \emph{conjugation rack} if it is a subset of a finite group $G$ that is closed under conjugation.

A \emph{subrack} of a rack $X$ is a subset of $X$ which is a rack. It is an easy fact that the set of all subracks of a finite rack $X$, denoted $\mathcal{R}(X)$, form a lattice under set-theoretical inclusion. It was proved in \cite[Theorem~1.1]{HSW19} that the subrack lattice $\mathcal{R}(G)$ of a finite group $G$ determines if the group $G$ is abelian, nilpotent, supersolvable, solvable or simple. In \cite{Kay21} Kayacan extended this result by proving the isomorphism type of the subrack lattice $\mathcal{R}(G)$ determines the nilpotence class of $G$.

A partially ordered set (abbreviated poset) $\mathcal{P}$ is called \emph{bounded} if it has a unique maximal element and a unique minimal element. For a bounded poset $\mathcal{P}$ with greatest element $\hat{1}$ and least element $\hat{0}$, we define the proper part of $\mathcal{P}$ as $\overline{\mathcal{P}} := \mathcal{P}\setminus \{\hat{1},\hat{0}\}$. To any bounded poset $\mathcal{P}$ one can associate its order complex $\Delta(\mathcal{P})$ which is the abstract simplicial complex of all linearly ordered subsets (chains) of $\overline{\mathcal{P}}$. In particular, the order complex $\Delta(\mathcal{R}(X))$ is the simplicial complex whose simplices are the chains of $\mathcal{R}(X)\setminus \{X,\emptyset\}$. Heckenberger, Shareshian, and Welker studied the homotopy properties of the subrack lattices of conjugation racks in the aforementioned paper. Another important result (see \cite[Proposition~1.3]{HSW19}) they obtained is that for a finite group $G$, the order complex $\Delta(\mathcal{R}(G))$ is homotopy equivalent to a $(c-2)$-sphere, where $c$ is the number of conjugacy classes of $G$.

Consider the order complex $\Delta(\mathcal{R}(X))$ of the subrack lattice of a finite rack $X$. The aim of this paper is to determine the homotopy type of $\Delta(\mathcal{R}(X))$ in some cases and, thereby, to generalize the known results on conjugation racks. Recall that for any $a\in X$, the map $\phi_a\colon X\to X;\,x\mapsto a\triangleright x$ is an automorphism of $X$. There is more to say. The map
\begin{align*}
  \Phi \colon X  & \longrightarrow\; \Phi(X) \\
  a & \longmapsto\; \phi_a
\end{align*}
defines a rack homomorphism when the image set $\Phi(X)$ is considered as a conjugation rack and the \emph{inner automorphism group} $\mathsf{Inn}(X) := \langle \Phi(X) \rangle$ is a normal subgroup of the full automorphism group of $X$ as the equality $\sigma\phi_a\sigma^{-1} = \phi_{\sigma(a)}$ holds for any automorphism $\sigma$ of $X$. Further, for a conjugation rack $X$ the map $\Phi$ extends to a group homomorphism from $\langle X \rangle$ to $\mathsf{Inn}(X)$ so that $\mathsf{Inn}(X) \cong \langle X \rangle / Z(\langle X \rangle)$ (see \cite[Proposition~2]{Kay20}). 

Let $G$ be a group acting on a set $X$. We denote the orbit of an element $a$ of $X$ under the action of $G$ with $\mathsf{Orb}_G(a)$. Notice that the set of $G$-orbits of $X$ form a partition of $X$. Let $S$ be a subset of $X$. We denote the stabilizer of $S$ under the action of $G$ with $\mathsf{Stab}_G(S)$ and if $S = \{a\}$ is a singleton we will just use $\mathsf{Stab}_G(a)$.

\begin{thm}\label{thm:thmA}
  Let $X$ be a finite rack having the following property. The poset $\overline{\mathcal{R}(X)}$ is non-empty and the equality $\mathsf{Stab}_{\mathsf{Inn}(X)}(M) = \mathsf{Inn}(X)$ holds for every maximal subrack $M$ of $X$. Then $\Delta(\mathcal{R}(X))$ is homotopy equivalent to a $(m-2)$-sphere, where $m$ is the number of $\mathsf{Inn}(X)$-orbits in $X$.
\end{thm}

Recall that for a group rack the homotopy type of the associated order complex is a $(c-2)$-sphere, where $c$ is the number of conjugacy classes of the group. Another similar result can be obtained for the conjugacy class racks in nilpotent groups. Actually, those two results can be deduced as corollaries of Theorem~\ref{thm:thmA} (see Corollary~\ref{cor:groups} and Corollary~\ref{cor:nil}).

For a non-empty rack $X$, we say $X$ is \emph{decomposable} if it admits a partition $X = Y\sqcup Z$ such that $Y$ and $Z$ are non-empty subracks; otherwise, we say $X$ is \emph{indecomposable}. Equivalently, $X$ is indecomposable if and only if $X = \mathsf{Orb}_{\mathsf{Inn}(X)}(a)$ for every $a\in X$. Notice that any rack $X$ having the property stated in Theorem~\ref{thm:thmA} is necessarily a decomposable rack.

For an integer $n$, the partition lattice $\Pi_n$ is the set of all partitions of $[n]:=\{1,2,\dots,n\}$ ordered by refinement. For $2\leq k\leq n$, the $k$-equal partition lattice $\Pi_{n,k}$ is the lattice of all partitions $B_1| \dots | B_r$ of $[n]$ such that, for all $1\leq i\leq r$, either $|B_i| = 1$ or $|B_i|\geq k$. In \cite[Example~2.3]{HSW19} we see that if $T$ is the rack of transpositions of the symmetric group $S_n$, then $\mathcal{R}(T)$ is isomorphic to the partition lattice $\Pi_n$. And in \cite[Proposition~2.11]{HSW19} it was proved that if $P$ is the rack of all $p$-cycles in the alternating group $A_n$, with $p < n-2$ being an odd prime number, then the order complexes of $\mathcal{R}(P)$ and $\Pi_{n,p}$ are homotopy equivalent. Clearly, both $T$ and $P$ are indecomposable racks.

Given an integer $m$ satisfying $mk \leq n$ we define the subposet $\Pi_{n,k}^m$ of $\Pi_{n,k}$ in the following way. An element of $\Pi_{n,k}$ is an element of $\Pi_{n,k}^m$ if and only if either it is the partition consisting of singletons or the number of elements that are lying in non-singletons blocks is at least $mk$. With this notation we can state our second main result. 

\begin{thm}\label{thm:thmB}
  Let $n = sp + t$ be an integer with $p$ being a prime number. For the possible values of $s$, $p$, and $t$ suppose one of the following cases occurs:
  \begin{tasks}[label=(\alph*),label-width=4ex](2)
 \task $s = 1$, $p\geq 2$ and $t > 2$
 \task $s = 2$, $p\geq 5$ and $t > 2$
 \task $s = 3$, $p\geq 5$ and $t > 3$
 \task $s = 4$, $p = 5$ and $t > 5$
 \task $s = 4$, $p\geq 7$ and $t > 4$
 \task $s = 5$, $p\geq 7$ and $t > 6$
 \task $s = 6$, $p\geq 11$ and $t > 6$
 \task $s = 7$, $p\geq 11$ and $t > 8$
 \task $s \geq 8$, $p\geq 2s - 1$ and $t > 4s - 4$
 \end{tasks}
 Let $C$ be the conjugacy class of elements in the symmetric group $S_n$ whose cycle type is $(p^s,1^t)$. Then the order complexes of $\mathcal{R}(C)$ and $\Pi_{n,p}^s$ are homotopy equivalent. 
\end{thm}

Notice that $\Pi_{n,2}$ is same with $\Pi_{n}$ and also $\Pi_{n,k}^1$ is same with $\Pi_{n,k}$. Therefore, in the above Theorem the case (a) corresponds to the rack $T$ of all transpositions in $S_n$ when $p = 2$ and to the rack $P$ of all $p$-cycles in $A_n$ when $p$ is an odd prime number. As a result of Theorem~\ref{thm:thmB} we can determine the homotopy types of the order complexes of subrack lattices for more general classes of conjugacy class racks (see Corollary~\ref{cor:kozlov}).

\section{Decomposable racks}\label{sec:dec}

Throughout the paper we shall introduce a couple of order-preserving poset maps (poset homomorphisms) that are useful to define homotopy equivalences between the corresponding order complexes.

\begin{lem}[see {\cite[Corollary~10.12]{Bjo95}}]\label{lem:poset}
Let $\mathcal{P}$ be a bounded poset and let $\epsilon\colon\mathcal{P}\to\mathcal{P}$ be a poset map whose restriction onto $\overline{\mathcal{P}}$ is an endomorphism of $\overline{\mathcal{P}}$. If $x \leq \epsilon(x)$ for all $x\in \mathcal{P}$, then $\epsilon$ induces homotopy equivalence between $\Delta(\mathcal{P})$ and $\Delta(\epsilon(\mathcal{P}))$. If, further, $\epsilon^2(x) = \epsilon(x)$ for all $x\in\mathcal{P}$, then $\Delta(\epsilon(\mathcal{P}))$ is a strong deformation retraction of $\Delta(\mathcal{P})$.
\end{lem}

A \emph{closure operator} $\epsilon$ on a bounded poset $\mathcal{P}$ is an endomorphism of $\mathcal{P}$ satisfying the following requirements:
\begin{enumerate}[(i)]
\item $\epsilon(\overline{\mathcal{P}}) \subseteq \overline{\mathcal{P}}$
\item $x \leq \epsilon(x)$ for all $x\in \mathcal{P}$
\item $\epsilon^2(x) = \epsilon(x)$ for all $x\in\mathcal{P}$
\end{enumerate}
Let $X$ be a finite rack and consider the mapping $\epsilon_X\colon \mathcal{R}(X)\to \mathcal{R}(X)$ taking a subrack $S$ to $\overline{S}$, where $\overline{S}$ is the union of $\mathsf{Inn}(X)$-orbits of the elements in $S$. Clearly, the mapping $\epsilon_X$ is order-preserving and it satisfies the above Conditions~(ii) and (iii) as well. However, Condition~(i) may not hold. Certainly it does not hold if $X$ is indecomposable and its cardinality is greater than 1.

\begin{lem}\label{lem:sub}
  A subset $S$ of a rack $X$ is a subrack if and only if
  $$ \langle \Phi(S) \rangle \leq \mathsf{Stab}_{\mathsf{Inn}(X)}(S). $$
\end{lem}

\begin{proof}
  If $S$ is a subrack of $X$, then for every $s\in S$ the map $\phi_s$ restricted to $S$ is a bijection of $S$ meaning every generator element of $\langle \Phi(S) \rangle$ lies in the stabilizer of $S$.

  Conversely, if $\langle \Phi(S) \rangle$ is contained in $\mathsf{Stab}_{\mathsf{Inn}(X)}(S)$, then $\phi_s^{\pm 1}(t)$ is an element of $S$ for every $s,t\in S$.
\end{proof}

\begin{prop}\label{prop:max}
  Let $X$ be a rack and $S$ be a subrack of $X$. Consider the following statements.
  \begin{enumerate}[(i)]
  \item The subrack $S$ is proper in $X$ and for each $\mathsf{Inn}(X)$-orbit $O$ the intersection $S\cap O$ is non-empty.
  \item The subgroup $\mathsf{Stab}_{\mathsf{Inn}(X)}(S)$ is proper in $\mathsf{Inn}(X)$ and there is no proper normal subgroup of $\mathsf{Inn}(X)$ containing $\mathsf{Stab}_{\mathsf{Inn}(X)}(S)$.
  \end{enumerate}
  Then the first statement implies the latter. Moreover, if $S$ is a maximal subrack of $X$ and $\mathsf{Stab}_{\mathsf{Inn}(X)}(S) \neq \mathsf{Inn}(X)$, then for each $\mathsf{Inn}(X)$-orbit $O$ the intersection $S\cap O$ is non-empty.
\end{prop}

\begin{proof}
  Let $X = \bigsqcup O_i$ be the partition of $X$ into $\mathsf{Inn}(X)$-orbits. By assumption, there exists $y\in X\setminus S$. Let $y\in O_i$ for some $\mathsf{Inn}(X)$-orbit $O_i$ and let $x\in S\cap O_i$. Thus, there exist a $\sigma\in\mathsf{Inn}(X)$ so that $\sigma(x) = y$. Clearly $\sigma\notin \mathsf{Stab}_{\mathsf{Inn}(X)}(S)$; hence, $\mathsf{Stab}_{\mathsf{Inn}(X)}(S)$ is a proper subgroup of $\mathsf{Inn}(X)$. Next, let $N$ be a normal subgroup of $\mathsf{Inn}(X)$ containing $\mathsf{Stab}_{\mathsf{Inn}(X)}(S)$. By Lemma~\ref{lem:sub}, the inner automorphism $\phi_x$ lies in $N$. Since $N$ is a normal subgroup of $\mathsf{Inn}(X)$ and $\phi_y = \phi_{\sigma(x)} = \sigma\phi_x\sigma^{-1}$, we see that $\phi_y$ also lies in $N$. However, this observation is valid for every $y\in X\setminus S$ implying $N = \mathsf{Inn}(X)$, since $N$ contains every inner automorphism.
  
  To prove the last assertion of the Proposition, let us suppose $S$ is a maximal subrack of $X$ and $\mathsf{Stab}_{\mathsf{Inn}(X)}(S) \neq \mathsf{Inn}(X)$. Observe that, if $S\cap O$ is empty for some $\mathsf{Inn}(X)$-orbit $O$, then $S = X\setminus O$ since $S$ is a maximal subrack of $X$. However, $\mathsf{Stab}_{\mathsf{Inn}(X)}(S) = \mathsf{Inn}(X)$ in that case which yields a contradiction.
\end{proof}

\begin{proof}[Proof of Theorem~\ref{thm:thmA}]
  By assumption $\mathsf{Stab}_{\mathsf{Inn}(X)}(M) = \mathsf{Inn}(X)$ for any maximal subrack $M$ of $X$; therefore, the maximal subracks are exactly the unions of all but one $\mathsf{Inn}(X)$-orbit in $X$ by Proposition~\ref{prop:max}. Consider the endomorphism $\epsilon_X$ of $\mathcal{R}(X)$ taking a subrack $S$ of $X$ to $\overline{S}$, where $\overline{S}$ is the union of $\mathsf{Inn}(X)$-orbits of elements in $S$. Since $\epsilon_X(S)\neq X$ unless $S=X$, we see that $\epsilon_X$ is a closure operator on $\mathcal{R}(X)$ and this implies $\mathcal{R}(X)$ and $\epsilon_X(\mathcal{R}(X))$ are of the same homotopy type as order complexes by Lemma~\ref{lem:poset}. Since $\epsilon_X(\mathcal{R}(X))$ can be identified with $2^{[m]}$, where $m$ is the number of $\mathsf{Inn}(X)$-orbits in $X$ and since $\Delta(2^{[m]})$ is the triangulation of a $(m-2)$-sphere, we are done. 
\end{proof}

\begin{cor}[see {\cite[Proposition~1.3]{HSW19}}]\label{cor:groups}
  Let $G$ be a finite group rack whose cardinality is greater than one. Then $\Delta(\mathcal{R}(G))$ is a $(c-2)$-sphere, where $c$ is the number of  conjugacy classes of $G$.
\end{cor}

\begin{proof}
  Observe that the conjugacy classes of the group $G$ are exactly the $\mathsf{Inn}(G)$-orbits of the elements of the group rack $G$. It is a standard fact in finite group theory that any subset of $G$ containing a representative element from each conjugacy class generates the whole group $G$. Therefore, the maximal subracks of the group rack $G$ are exactly the unions of all but one conjugacy class of $G$ and Theorem~\ref{thm:thmA} applies.
\end{proof}

Another similar result holds for the conjugacy classes of nilpotent groups. First we shall prove an auxiliary result.

\begin{lem}\label{lem:p-not-conn}
  In a finite nilpotent group, a conjugacy class of a non-central element is decomposable as a rack.
  
\end{lem}

\begin{proof}
  Let $G$ be a counter example, i.e. $G$ is a nilpotent group and there exist a conjugacy class $C$ of $G$ such that $|C|>1$ and $\mathsf{Inn}(C)$ acts on $C$ transitively. Since $\mathsf{Inn}(C)\cong \langle C\rangle / Z(\langle C\rangle)$, the action of $\langle C\rangle$ on $C$ by conjugation is transitive as well.
  
\medskip
\noindent
\emph{Claim.} In a finite nilpotent group which is not cyclic of prime power order, subgroups generated by the conjugacy classes are proper.

\medskip
This follows from the facts that any element is contained in a maximal subgroup (except in a cyclic group of prime power order) and maximal subgroups are normal in nilpotent groups.

\medskip
Now, by the Claim, $H:=\langle C\rangle$ is a proper subgroup of $G$. Clearly, $C$ is a conjugacy class of $H$. For $H$ being a nilpotent group, the above Claim implies $\langle C\rangle < H$ which is a contradiction.
\end{proof}

\begin{cor}\label{cor:nil}
  Let $G$ be a finite nilpotent group and $C$ be a conjugacy class of a non-central element in $G$. Then $\Delta(\mathcal{R}(C))$ is a $(m-2)$-sphere, where $m$ is the number of maximal elements of $\mathcal{R}(C)$. 
\end{cor}

\begin{proof}
  We want to show that the assumption of Theorem~\ref{thm:thmA} applies for the conjugacy class rack $C$, i.e., we want to prove that the maximal subracks of $C$ are the unions of all but one $H$-orbit of elements of $C$, where $H:=\langle C\rangle$.

  Observe that $C=\bigsqcup_{i=1}^mC_i$ is the disjoint union of some $H$-orbits  $C_i,\; 1\leq i\leq m$. By Lemma~\ref{lem:p-not-conn}, the number $m$ of those orbits is greater than one. Suppose contrarily there exist a maximal subrack $M$ of $C$ which is not the union of $H$-orbits of $C$ except one. If $M\cap C_i=\emptyset$ for some $1\leq i\leq m$, then $C\setminus C_i$ is a proper subrack of $C$ containing $M$ contrary to our assumption. Therefore $M\cap C_i\neq\emptyset$ for any $1\leq i\leq m$. Let $A$ be a $H$-orbit which is not contained by $M$.

  Since any $H$-orbit is actually a conjugacy class of $H$, if a normal subgroup $N$ of $H$ contains an element of $C_i$ for some $1\leq i\leq m$, then it contains all the elements of $C_i$. Consider the normalizer $N_H(M)$. If $N_H(M)=H$, then $M\cap A$ would be a $H$-orbit contrary to our assumption. Thus $N_H(M)$ is a proper subgroup of $H$. Next suppose that $N_H(M)$ is a normal subgroup of $H$. Then $C\subset N_H(M)$ as $M\cap C_i\neq\emptyset$ for any $1\leq i\leq m$ implying $N_H(M)\geq \langle C\rangle=H$ which is a contradiction. Hence $N_H(M)$ is a proper non-normal subgroup of $H$. As $H$ is a nilpotent group, maximal subgroups of $H$ are normal; so there exist a normal subgroup $K$ of $H$ containing $N_H(M)$ properly. However, this yields again a contradiction as $K\cap C_i\neq\emptyset$ for any $1\leq i\leq m$ but $K$ is a proper normal subgroup of $H$ by assumption.  
\end{proof}

In the statement of Theorem~\ref{thm:thmA} we made the assumption that ``the equality $\mathsf{Inn}(X) = \mathsf{Stab}_{\mathsf{Inn}(X)}(M)$ holds for every maximal subrack $M$ of $X$''. Alternatively, using the first assertion of Proposition~\ref{prop:max} in the contrapositive form, we may rewrite Theorem~\ref{thm:thmA} with the assumption that ``for a subrack $S$ of $X$ either $\mathsf{Stab}_{\mathsf{Inn}(X)}(S) = \mathsf{Inn}(X)$ or $\mathsf{Stab}_{\mathsf{Inn}(X)}(S)$ is contained by a proper normal subgroup of $\mathsf{Inn}(X)$''. Deducing Corollary~\ref{cor:nil} would be easier with this alternative reformulation of Theorem~\ref{thm:thmA}.

\section{Indecomposable racks}
\label{sec:conn}

A \emph{Galois connection} between two posets $\mathcal{P}$ and $\mathcal{Q}$ is a pair of order-preserving mappings $\alpha\colon \mathcal{P}\to \mathcal{Q}$ and $\beta\colon \mathcal{Q}\to \mathcal{P}$ having the following properties:
\begin{enumerate}[(i)]
\item $x\leq \beta\alpha(x)$ for all $x\in \mathcal{P}$
\item $\alpha\beta(y)\leq y$ for all $y\in \mathcal{Q}$
\end{enumerate}
Suppose the pair $\alpha\colon \mathcal{P}\to \mathcal{Q}$ and $\beta\colon \mathcal{Q}\to \mathcal{P}$ is a Galois connection between $\mathcal{P}$ and $\mathcal{Q}$. If $\mathcal{P}$ and $\mathcal{Q}$ are bounded posets and if $\alpha(\overline{\mathcal{P}}) \subseteq \overline{\mathcal{Q}}$ and $\beta(\overline{\mathcal{Q}}) \subseteq \overline{\mathcal{P}}$, we can mention about a Galois connection between $\overline{\mathcal{P}}$ and $\overline{\mathcal{Q}}$ which is the restriction of those maps onto the proper parts. Suppose $\mathcal{P}$ is a bounded poset with greatest element $\hat{1}$ and suppose $\beta\alpha(x) = \hat{1}$ if and only if $x = \hat{1}$. Then the mapping $\gamma\colon \mathcal{P}\to \mathcal{P};\,x\mapsto \beta\alpha(x)$ would be a closure operator \cite{Rot64}.

For a finite group $G$, we denote by $\mathcal{L}(G)$ the subgroup lattice of $G$. Let $X$ be a finite rack. Associated with $X$ we define two functions. The mapping $\alpha_X\colon\mathcal{R}(X)\to\mathcal{L}(\mathsf{Inn}(X))$ takes a non-empty subrack $S$ to $\langle\Phi(S)\rangle$ and the image of the empty rack $\emptyset$ is the trivial subgroup $1$. And the mapping $\beta_X\colon\mathcal{L}(\mathsf{Inn}(X))\to\mathcal{R}(X)$ takes a subgroup $H$ to the subrack $S$, where $S$ is the largest subrack satisfying $\langle\Phi(S)\rangle\leq H$. Clearly, $\alpha_X$ and $\beta_X$ are order-preserving mappings. Notice that when $X$ is indecomposable the inequality $\alpha_X(S)\neq \mathsf{Inn}(X)$ holds for every proper subrack $S$ of $X$. Therefore, for an indecomposable rack $X$ the mapping $\gamma_X\colon \mathcal{R}(X)\to \mathcal{R}(X);\,S\mapsto \beta_X\alpha_X(S)$ is a closure operator.

\begin{rem}
 Let $C$ be a conjugacy class of a group $G$ such that $G = \langle C \rangle$ and $Z(G) = 1$. By \cite[Proposition~2]{Kay20} the conjugation rack $C$ is indecomposable and faithful. Further, the inner automorphism group $\mathsf{Inn}(C)$ is isomorphic with $G$ in a natural way so that, by replacing $\mathsf{Inn}(C)$ with $G$, we can define the mappings $\alpha'_C\colon\mathcal{R}(C)\to\mathcal{L}(G)$ and $\beta'_C\colon\mathcal{L}(G)\to\mathcal{R}(C)$ in the obvious way. Let $S$ be a subrack of $C$ and let $\gamma'_C\colon\mathcal{R}(C)\to\mathcal{R}(C)$ be the mapping taking $S$ to $\beta'_C\alpha'_C(S)$. By the above discussion it should be clear that $\gamma'_C(S) = C\cap\langle S\rangle$. 
\end{rem}

It is difficult to mention about a Galois connection between the proper parts of $\mathcal{R}(X)$ and $\mathcal{L}(\mathsf{Inn}(X))$ as $\beta_X(H)$ may be the empty rack $\emptyset$ for a non-trivial subgroup $H$ of $\mathsf{Inn}(X)$. In the following example it is easy to describe a subposet of the subgroup lattice so that the restriction of the mappings $\alpha_X$ and $\beta_X$ yields a Galois connection.

\begin{example}
 Let $n$ be an integer and $p$ be an odd prime such that $n/2 < p < n-1$. Let $C$ be the conjugacy class of $p$-cycles in the alternating group $A_n$. Observe that $C$ is an indecomposable faithful rack and we can identify $\mathsf{Inn}(C)$ with $A_n$ as in the previous remark. Let $\mathcal{L}_p(A_n)$ be the subposet of $\mathcal{L}(A_n)$ consisting of subgroups of $A_n$ whose orders are  divisible by $p$. Let the mappings $\alpha''_C\colon\mathcal{R}(C)\to\mathcal{L}_p(A_n)$ and $\beta''_C\colon\mathcal{L}_p(A_n)\to\mathcal{R}(C)$ be defined in the obvious way. Then the pair $\alpha''_C$ and $\beta''_C$ form a Galois connection between $\mathcal{R}(C)$ and $\mathcal{L}_p(A_n)$.
\end{example}

Recall that the mapping $\gamma_X\colon \mathcal{R}(X)\to \mathcal{R}(X)$ is an endomorphism of $\mathcal{R}(X)$ which is a closure operator if $X$ is indecomposable. Orbit decompositions are also useful to define a poset endomorphism. For a subrack $S$ of $X$ the \emph{orbit structure} of $S$ is the partition of $X$ into $\langle \Phi(S) \rangle$-orbits. Consider the mapping $\eta_X\colon \mathcal{R}(X)\to \mathcal{R}(X)$ taking a subrack $S$ into $\widetilde{S}$, where $\widetilde{S}$ is the largest subrack of $X$ having the same orbit structure with $S$.

\begin{lem}
 Let $X$ be a finite rack. Then the mapping $\eta_X\colon \mathcal{R}(X)\to \mathcal{R}(X)$ described above is a well-defined poset endomorphism. Moreover, it is a closure operator on $\mathcal{R}(X)$ when $X$ is indecomposable.
\end{lem}

\begin{proof}
 Let $S$ be a subrack of $X$. For an element $a$ of $X$, the orbits of $\langle \Phi(S)\rangle$ is preserved by the automorphism $\phi_a$ if and only if the orbit structure of $\{a\}$ is same with or more refined than the orbit structure of $S$. And the set of all elements of $X$ having this property constitutes a subrack. This shows that $\eta_X$ is well-defined. Also, it is easy to see that $\eta_X$ is an order-preserving map. Finally, if $X$ is an indecomposable rack, then $\widetilde{S} = X$ implies $S = X$ as $X$ is the single $\langle \Phi(S)\rangle$-orbit in that case. 
\end{proof}

Suppose that $X$ is an indecomposable finite rack. Notice that each element of $\eta_X(\mathcal{R}(X))$ can be uniquely identified with a partition of $X$. That means $\eta_X(\mathcal{R}(X))$ can be considered as a subposet of $\Pi_{|X|}$.

\begin{rem}
 Let $C$ be a conjugacy class of a transitive permutation group $G$ of degree $n$ such that $G = \langle C \rangle$ and $Z(G) = 1$. As was previously remarked the conjugation rack $C$ is an indecomposable, faithful rack and the inner automorphism group $\mathsf{Inn}(C)$ is isomorphic with $G$ in a natural way. Let $\eta'_C\colon\mathcal{R}(C)\to\mathcal{R}(C)$ be the mapping taking a subrack $S$ to $\widehat{S}$, where $\widehat{S}$ is the largest subrack of $C$ determining the same partition of $[n]$ with the partition given by the $\langle S\rangle$-orbits of the elements of $[n]$. Notice that $\eta'_C$ is a well-defined poset endomorphism. However, it may be the case that $\widehat{S} = C$ for a proper subrack $S$ of $C$. 
 
 
\end{rem}

\begin{lem}[see {\cite[Lemma~2.10]{HSW19}}]\label{lem:hat}
 Let $C$ be a conjugacy class rack described as in Theorem~\ref{thm:thmB}. Then, for a subrack $S$ of $C$, we have the equality $\widehat{S} = C$ if and only if $S=C$.
\end{lem} 

\begin{proof}
 We call the symmetric group $S_n$ or the alternating group $A_n$ the \emph{giant} permutation groups on $[n]$. Since the inner automorphism group of the rack $C$ is isomorphic to a giant, we see that the rack $C$ is indecomposable and $\langle C\rangle$ acts transitively on $[n]$.
 
 We want to show that if $\langle S\rangle$ acts on $[n]$ transitively, then $\langle S\rangle$ is a giant. That would imply the equality $S = C$ holds. By \cite[Theorem~13.10]{Wie64} $\langle S\rangle$ is a giant when $\langle S\rangle$ acts primitively on $[n]$. In the rest of the proof we show that if $\langle S\rangle$ acts transitively on $[n]$ then its action is primitive as well.
 
 Suppose contrarily there is a block system $B_1| \dots | B_q$ for the action of $\langle S\rangle$ on $[n]$. Since the action is transitive by assumption, the size of each block is same and it is at least $p$ as the cycle type of a generator element is $(p^s,1^t)$. Let $\pi$ be an element of $S$ which does not fix one of the blocks $B_i$, $1\leq i\leq q$. As the cycle type of $\pi$ is $(p^s,1^t)$, we see that $\pi$ must move at least $p$ blocks. However, that means the value of $s$ is at least $p$ which is not possible. 
\end{proof}

\begin{proof}[Proof of Theorem~\ref{thm:thmB}]
 By the previous Remark and Lemma~\ref{lem:hat} we know that the map $\eta'_C\colon\mathcal{R}(C)\to\mathcal{R}(C)$ is a closure operator. Moreover, each element of $\eta'_C(\mathcal{R}(C))$ can be uniquely identified with an element $\Pi_n$ in a natural way so that the identification $\iota\colon \eta'_C(\mathcal{R}(C))\to \Pi_n$  is an injective poset map. Therefore, to prove the Theorem it is enough to show that $\iota(\eta'_C(\mathcal{R}(C))) = \Pi_{n,p}^s$. It is clear that for any $S\in \eta'_C(\mathcal{R}(C))$ the partition $\iota(S)$ is an element of $\Pi_{n,p}^s$. Suppose $\pi\in \Pi_{n,p}^s$. Let $B_i$, $1\leq i\leq r$, be the blocks of $\pi$. Notice that the number of non-singleton blocks is at least $s$ as $\pi\in \Pi_{n,p}^s$. Let $S$ be the set of all permutations in $C$ that are fixing the blocks of $\pi$. Notice that $S\neq \emptyset$. Since the action of $S$ is transitive on each block $B_i$ and since $S$ is the largest subrack of $C$ with this property, we see that $\iota^{-1}(\pi) = S$. 
\end{proof}

A finite simplicial complex $\Delta$ is said to be \emph{shellable} if its maximal faces can be arranged in linear order $F_1,F_2,\dots,F_t$ in such a way that the subcomplex $\left(\bigcup_{i=1}^{k-1}F_i\right) \cap F_k$ is pure and ($\dim F_k - 1$)-dimensional for all $k = 2,3,\dots,t$. It is a standard fact that the homotopy type of a shellable complex is a wedge of spheres. Let $\mathcal{P}$ be a finite poset. There are various ``lexicographic'' labeling rules for the edge set of the Hasse diagram of $\mathcal{P}$ such as EL-labeling and EC-labeling (see \cite{BW96,Koz97}). The existence of any of those labelings implies the shellability of the order complex $\Delta(\mathcal{P})$. Accordingly, we say $\Delta(\mathcal{P})$ is \emph{EL-shellable} (respectively, \emph{EC-shellable}) if $\mathcal{P}$ admits an EL-labeling (respectively, EC-labeling).

\begin{cor}\label{cor:kozlov}
  Consider the possible cases in Theorem~\ref{thm:thmB}.
  \begin{itemize}
   \item Case (a): The order complex $\Delta(\mathcal{R}(C))$ is EL-shellable and it is homotopy equivalent to a wedge of spheres. 
   \begin{itemize}
    \item ($p > t$): $\Delta(\mathcal{R}(C))$ is homotopy equivalent to a wedge of ($t-1$)-spheres. 
    \item ($p \leq t$): There are spheres of several dimensions in the wedge.
   \end{itemize}
   \item Cases (b)--(i): 
   \begin{itemize}
    \item ($p > t$): The order complex $\Delta(\mathcal{R}(C))$ is EC-shellable and it is homotopy equivalent to a wedge of ($s+t-2$)-spheres.
   \end{itemize}
   \pagebreak
   \item Cases (c)--(i): 
   \begin{itemize}
    \item ($p \leq t$): The order complex $\Delta(\mathcal{R}(C))$ is not shellable.
   \end{itemize}
   
  \end{itemize}

\end{cor}

\begin{proof}
 By Theorem~\ref{thm:thmB} the order complexes of $\mathcal{R}(C)$ and $\Pi_{n,p}^s$ are homotopy equivalent. Then the proof follows from \cite[Theorem~5.3]{Koz97}, \cite[Theorem~5.9]{BW96} and \cite[Theorem~3.9]{Koz97}.
\end{proof}



\end{document}